# AHAB'S ARITHMETIC;
## OR, THE MATHEMATICS OF MOBY-DICK

Sarah Hart; Birkbeck College, University of London; s.hart@bbk.ac.uk

## Abstract

Herman Melville's novel *Moby-Dick* contains a surprising number of mathematical allusions. In this article we explore some of these, as well as discussing the questions that naturally follow: why did Melville choose to use so much mathematical imagery? How did Melville come to acquire the level of mathematical knowledge shown in the novel? And is it commensurate with the general level of mathematical literacy at that time?

## 1. Introduction

*Moby-Dick* is one of the most famous novels ever written. 'It is a surpassingly beautiful book,' wrote D.H. Lawrence in 1922. 'It is a great book, a very great book, the greatest book of the sea ever written. It moves awe in the soul.' Melville deploys a vast array of literary, cultural and religious references, symbolism and imagery.

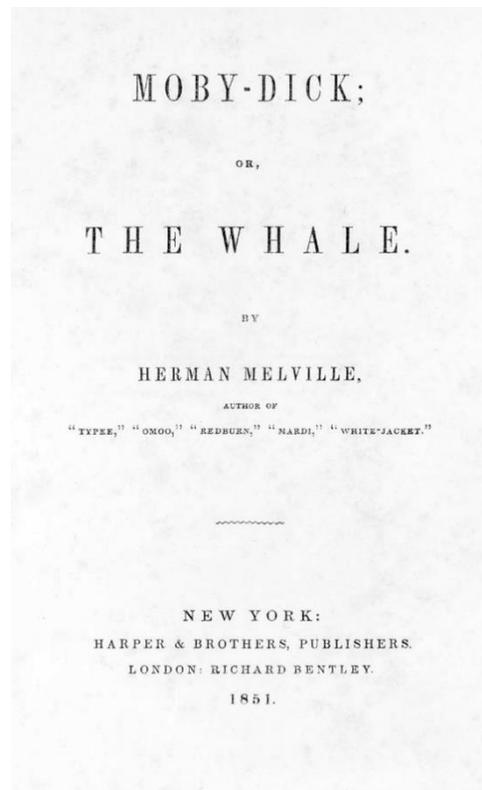

**Figure 1: Title page of the first edition of *Moby-Dick*, 1851 (Beinecke Library, Yale University)**

Any mathematician reading *Moby-Dick* would be struck, as I was, by the number of mathematical references, and evident mathematical knowledge, shown by that book (several examples will be given in this paper). But is the number of mathematical references really unusual? What were contemporary authors doing and saying? Was Melville especially mathematical? If so, why? I will argue that yes, Melville was indeed atypical in his level of mathematical knowledge, and moreover that he evidently enjoyed mathematics and was good at it. The main motivation for this article is to explore these questions, but it is hoped that the examples and quotations discussed could also be useful pedagogically, as a way for educators to contextualise some of the mathematics which is referenced. Finally, we hope that this article reinforces the idea of Melville as being another counterexample, were one needed, to the false 'two cultures' dichotomy.

Melville must be one of the most written-about authors, and *Moby-Dick* one of the most written-about books, of all time. There are many biographies, and numerous journal articles – indeed at least one dedicated journal, *Leviathan*, of Melville studies. Yet I have only been able to find two papers (Turpin, 2015; Farmer, 2016) among the hundreds written by Melville scholars that address this aspect of Melville's work – neither of them explicitly focuses on the mathematics of *Moby-Dick*. Turpin's essay looks at instances of Melville's allusions to mathematics in the context of his (Melville's) interest in Platonic idealism. Meredith Farmer's excellent article on Melville's education (Farmer, 2016) sheds light on the unusual opportunities Melville had for learning science and mathematics; the essay convincingly refutes the prevailing view among his biographers that Melville 'had no more instruction in science' 'than the average boy of his times'. Mathematical references can be found in many, perhaps all, of Melville's novels. In *Mardi*, for example, he has the character Babbalanja cry out '*Oh Man, Man, Man! Thou art harder to solve, than the Integral Calculus*'. Elsewhere another character, in frustration at Babbalanja's philosophising, retorts: '*Away with your logic and conic sections!*' *Mardi*, though, had such a terrible critical reception that Melville promised his publisher that his next book would contain 'no metaphysics, no conic-sections, nothing but cakes & ale' (see Turpin, 2015). This may have been true for that next book, *Redburn*, but by the time *Moby-Dick* came along, Melville had most assuredly reneged on his promise. *Moby-Dick* seems to me the most mathematical of all Melville's books, and it is therefore on *Moby-Dick* that this essay concentrates.

We start in Section 2 with a brief biography of Herman Melville; in Section 3 we give an overview of *Moby-Dick* and its critical reception. Section 4 looks at some of the ways mathematical imagery is deployed in the novel, as evidence that Melville was an unusually mathematically-minded writer. We then look in Section 5 in more detail at Melville's education, which, due to an exceptional teacher, was likely to have featured more interesting mathematics than the standard schooling of the time. Section 6 looks at the mathematicians and mathematics books that are mentioned in *Moby-Dick*.

Meredith Farmer writes (Farmer, 2016) that one of the goals of her essay is to 'be available as historical scaffolding for work on Melville's engagement with different branches of science, technology, engineering and mathematics'. It is my hope that the current article can form a small contribution to this work.

# 2. A brief biography

Herman Melville was a very private person – here was a man who would hang a towel over the doorknob of his study in order to obscure the keyhole in case anyone might be tempted to peep through it. Only around 300 of his letters survive – as Delbanco points out, this is tiny in comparison to the 12,000 surviving letters by Henry James (Delbanco, 2006). Add to this the fact that after some early success his work was badly received and he sank into obscurity and relative poverty even in his own lifetime, and the challenge for any biographer is magnified. I have drawn for this biographical sketch mainly on the books of Delbanco and Parker ((Delbanco, 2006; Parker, 1996).

Herman Melville was born in New York on 1st August 1819, so 2019 is his bicentenary. Herman was the third of eight children of Allan and Maria Melvill (the "e" was added by Maria after Allan's death). Maria came from the respected and well-off Gansevoort family; she was a very religious woman and the bible was read in the house every day. Allan was an import merchant, whose many business enterprises invariably ended up in failure and debt. He had to borrow repeatedly from his and Maria's families. By October 1830 Allan Melvill's financial affairs were in such a bad state that the family had to leave New York, leaving many debts behind them. They moved to Albany, where they were financially supported by Maria's brother Peter Gansevoort. By the time Allan died in 1832 there was no choice but for Herman to find work, and he got a job as an errand boy at the New York State bank in Albany. He only had a few more months of formal education after this, in 1835 and again in 1836. He moved from job to job, never settling for long. He had a short-lived teaching job, then in Autumn 1838 did some training in engineering at a local academy (Lansingburgh), with a view to getting a surveying job in the construction of the Erie Canal, but this never materialised.

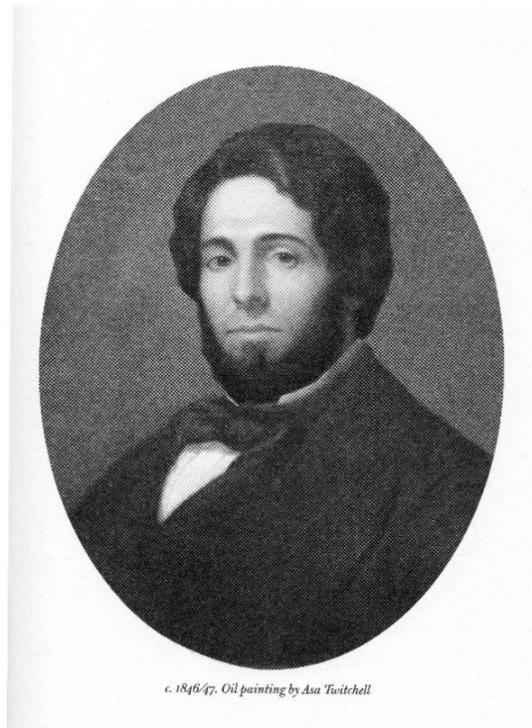

**Figure 2: Oil painting of Herman Melville, July 1846, by Asa Twitchell (Berkshire Athenaeum, Pittsfield, Mass.)**

In the Spring of 1839, frustrated with life in Albany, he followed the example of one of his cousins, and sailed as a cabin boy aboard a merchant ship travelling between New York and Liverpool. On his return home he tried teaching again, in the Greenbush village school in Albany, but this too didn't last. Eventually he reached what would be a life-changing decision to sign up, on January 3rd, 1841, as a deckhand on the whaling ship Acushnet. As Melville has Ishmael say, '*a whale-ship was my Yale College and my Harvard*' [Ch.24, p.122] [1].

---

[1] When quotations are referred to in the text of this article, we give a chapter and page number. Page numbers refer to the Penguin Classics edition (Melville, 2003). All quotations from the text, along with many other mathematical quotations – as complete a list as I can compile – are given in (Hart, 2019).

Melville spent almost four years at sea on several different ships, at one point spending a month in the South Pacific with a Polynesian tribe called the Typee. There was also time in Tahiti, Honolulu and elsewhere, finally returning to Boston on October 3rd 1844. On his return, after encouragement from friends and family, he began writing a fictionalised account of his experiences with the Typee tribe. The book, *Typee*, was submitted in Summer 1845 to Harper Brothers, who immediately rejected it on the grounds that it was 'impossible it could be true'. The manuscript was eventually accepted by Wiley and Putnam, and appeared in February 1846. He followed it up with *Omoo* (the title is the transliteration of a Polynesian word meaning 'wanderer').

*Omoo* appeared in Spring 1847. Melville was married in August 1847, to Elizabeth Shaw, and they moved to New York with Herman's younger brother Allan and his new wife. From 1847 to 1849, Herman was very productive, writing most of *Mardi*, all of *Redburn* and *White-Jacket*, and in the first weeks of 1850 the first few chapters of a new book that would eventually metamorphose into *Moby-Dick*. Seafaring tales were popular at that time, and Melville's books, especially the first two, were well received as part of that type. However with each succeeding book the reception was less positive, and indeed Melville himself was rather scathing about his own work. About *Redburn* he said 'I, the author, know [it] to be trash, and wrote it to buy some tobacco with'. However, this was a fertile time creatively for Melville. As well as writing three novels, Melville was reading voraciously: 'I have swam through libraries', he said. He read Virgil, Milton, Shakespeare, Voltaire, Mary Shelley, Dante, Schiller, Thackeray and many others. It was also around this time that he met and became close friends with the author Nathaniel Hawthorne, who supplied one of the few descriptive remarks we have about Melville as a man: he was apparently, though a gentleman, 'a little heterodox in the matter of clean linen'.

Melville moved out of New York in the summer of 1850 and completed *Moby-Dick* over the next year. It was published on November 14th 1851. Reviews were, to say the least, mixed. It was completely different from, and hugely more ambitious than, anything he had produced before. Readers hoping for a standard seafaring yarn would have been sorely disappointed. It almost sunk without trace in his lifetime – his lifetime earnings from it were a grand total of $556.37.

Things got worse: his next novel *Pierre* was universally panned as was much of everything else he wrote subsequently. He wrote almost no prose after 1853, the memorable exceptions being the compelling but strange novella *Bartleby, the Scrivener* (1853) and the late work *Billy Budd* (1889), which remained unpublished until 1924. He spent most of the last two decades of his life working for the US customs service. He died in on September 28, 1891, largely forgotten.

# 3. *Moby-Dick* and its critical reception

The plot of *Moby-Dick* is well known, but we give a brief summary here of its main points and the key characters. The narrator ('Call me Ishmael') is employed as a deck-hand on a whaling-ship, the *Pequod*, with its captain, Ahab, first mate Starbuck (latterly of coffee-shop fame) and second mate Stubb. It becomes gradually clear that for Ahab the purpose of this voyage is revenge: he is hell-bent on hunting and killing Moby Dick, the great white whale who was the cause of his lost leg. This fixation gradually takes over the whole mission; Ahab becomes increasingly obsessive and irrational, ultimately sealing his own fate and that of his crew with him. Ishmael alone lives to tell the tale. Right from the start it is clear that this is no ordinary adventure narrative. The first chapters consist of 'excerpts' describing mentions of whales and whaling from a vast array of sources from the Bible, to Shakespeare, to books of natural history. Ishmael and others frequently break off into philosophical musings; there are numerous chapters on the anatomy and physiology of whales, and much else besides. The book is dense with meaning and metaphor, as Ishmael says it must be when the subject, Leviathan itself, is so huge. *'Give me a condor's quill!'* says he, *'Give me Vesuvius' crater for an inkstand!'* [Ch.104, p.497]

Literary criticism of *Moby-Dick* has a 170-year history. Nick Selby discusses this history in his critical guide (Selby, 1998), which includes lengthy extracts from several reviews. Early reviews were mixed. The review in Harper's New Monthly Magazine (Selby 1998; 23) stated that *Moby-Dick* 'in point of richness and variety of incident, originality of conception, and splendour of description, surpasses any of the former productions of this highly successful author. […] the genius of the author for moral analysis is scarcely surpassed by his wizard power of description.' By contrast the October 1851 review in the London Athenaeum (Selby 1998; 25) was damning: 'Mr Melville has to thank himself only if his horrors and his heroics are flung aside by the general reader, as so much trash belonging to the worst school of Bedlam literature'. Harrison Ainsworth, writing in 1853, says, '[t]he style is maniacal – mad as a March hare – mowing, gibbering, screaming, like an incurable Bedlamite, reckless of keeper or straight-waistcoat'. Rehabilitation for the novel would take 70 years, when it began to be championed by authors such as D.H. Lawrence (quoted earlier) and E.M. Forster, who described it as 'immensely important' (Selby, 1998; 42).

*Moby-Dick* was vigorously championed in the 20th century as an emblematic work of the 'American Renaissance' by scholars in the new field of American Studies. To contextualise this slightly, one of the challenges for the development of a truly 'American' literary tradition in the 19th century had been that international copyright laws were not yet established. This meant that American publishers could print editions of works by British authors without paying royalties and could therefore make much bigger profits from them. (This is partly why Dickens was such a sensation in America – he complained that he hadn't made a penny from the huge book sales in the US.) Between 1820 and 1830, only around 100 novels by American writers were published. This number roughly trebled in the next decade, and by the 1840s there were nearer a thousand American novels. Melville can thus be said to have been part of the first wave of American writers of fiction. His work, and *Moby-Dick* in particular, became, from the 1920s, key components in the newly established narrative of American literary history. Since that time, *Moby-Dick* has been viewed as one of the great works of literature, whose 'greatness is unlike that of any other book', as Leo Bersani writes (Selby, 1998; 134).

The novel has been interpreted in many ways, over a dozen of which appear in Slade's book on the symbolism of *Moby-Dick* (Slade, 1998), who remarks that *Moby-Dick* 'has enjoyed voluminous and penetrating criticism for a century'. The white whale has been interpreted, among other things, as a symbol for Milton's Satan, the Old Testament God, Death, evil, terror, good, fate, man's interpretation of the meaning of life, and resurrection. Sometimes it is even thought of as a whale. Ahab, meanwhile, is Faust, or abolitionist Garrison, or anti-abolitionist John Calhoun (Selby, 1998; 141), or, for D.H. Lawrence, 'the last phallic being of the white man'.

The book as a whole has been read as a study of monomania, of hubris, or even, in one memorable article 'the psychological and philosophical ramifications of contemporary theories of digestion' (Doty, 2017). Slade (Slade, 1998; 38) refers to an essay by the psychoanalyst Henry A. Murray in which he writes, '[s]ince Ahab has been proclaimed the "Captain of the Id", the simplest psychological formula for Melville's dramatic epic is this: an insurgent Id in mortal conflict with an oppressive cultural superego.' (A lengthy extract from this essay is given in (Selby, 1998)).

# 4. Mathematical Imagery and Symbolism in *Moby-Dick*

'Critics feel that the nature of one's experiences will determine the meaning he will see in the whale,' says Slade (Slade, 1998; 27). Thus, as a mathematician perhaps it was inevitable that I would spot mathematical allusions. In this section I will illustrate with selected examples that the choice of metaphors and imagery used by Melville shows the pleasure he seems to take in mathematical ideas. I will also look at the symbolic use of calculation and measurement in the novel as attempts to control the natural world.

## *Mathematical Metaphors*

Melville seems to enjoy thinking about how large numbers or quantities are expressed. When describing a whale's carcass being jettisoned, he writes: '*The vast white headless phantom floats further and further from the ship, and every rod that it so floats, what seem square roods of sharks and cubic roods of fowls, augment the murderous din*' [Ch. 69, p.336]. The very mast and contents of the ship are used later to represent the vast age of a crew member, Fedallah. '"*How old do you suppose Fedallah is, Stubb?" "Do you see that mainmast there?" pointing to the ship; "well, that's the figure one; now take all the hoops in the Pequod's hold, and string along in a row with that mast, for oughts, do you see; well, that wouldn't begin to be Fedallah's age. Nor all the coopers in creation couldn't show hoops enough to make oughts enough*"' [Ch.73, p.356].

Mathematical ideas creep into descriptions throughout the book. Right at the start of the story, Ishmael describes the parsimony of the landlord at the inn where he is staying, thus: '*Abominable are the tumblers into which he pours his poison. Though true cylinders without—within, the villainous green goggling glasses deceitfully tapered downwards to a cheating bottom. Parallel meridians rudely pecked into the glass, surround these footpads' goblets*' [Ch.3, p.15]. Towards the end of the novel, Ahab praises a loyal cabin boy with a geometrical metaphor: '*True art thou, lad, as the circumference to its centre*' [Ch.129, p.581].

In the chapters describing the physiology of whales, Ishmael describes the sperm whale as having a 'pervading dignity' because of the 'mathematical symmetry' of its head. Symmetry is thus being used to symbolise virtue. Detail is given on the shape of the whale's head, including what Melville claims in a footnote is a previously undefined mathematical term.

*'Regarding the Sperm Whale's head as a solid oblong, you may, on an inclined plane, sideways divide it into two quoins\* whereof the lower is the bony structure, forming the cranium and jaws, and the upper an unctuous mass wholly free from bones;*
*\* Quoin is not a Euclidean term. It belongs to the pure nautical mathematics. I know not that it has been defined before. A quoin is a solid which differs from a wedge in having its sharp end formed by the steep inclination of one side, instead of the mutual tapering of both sides.'* [Ch.77, p.371]

Euclid is also referenced in discussing how whales process images.

*'How is it, then, with the whale? True, both his eyes, in themselves, must simultaneously act; but is his brain so much more comprehensive, combining, and subtle than man's, that he can at the same moment of time attentively examine two distinct prospects, one on one side of him, and the other in an exactly opposite direction? If he can, then is it as marvellous a thing in him, as if a man were able simultaneously to go through the demonstrations of two distinct problems in Euclid. Nor, strictly investigated, is there any incongruity in this comparison.'* [Ch.74, p.361]

While gazing out to sea, Ishmael's contemplations often have a mathematical flavour. 'There you stand, lost in the infinite series of the sea, with nothing ruffled but the whales' [Ch.35, p.169]. He writes about watching a whale swimming far off in the ocean:

*'Even if not the slightest other part of the creature be visible, this isolated fin will, at times, be seen plainly projecting from the surface. When the sea is moderately calm, and slightly marked with spherical ripples, and this gnomon-like fin stands up and casts shadows upon the wrinkled surface, it may well be supposed that the watery circle surrounding it somewhat resembles a dial, with its style and wavy hour-lines graved on it. On that Ahaz-dial the shadow often goes back.'* [Ch.32, p.151]

This lovely 'fin as gnomon' metaphor is made even better when we discover that the 'Ahaz-dial' refers to what is believed to be the earliest written reference to sundials, in the book of Isaiah. '*Behold, I will bring again the shadow of the degrees, which is gone down in the sun dial of Ahaz, ten degrees backward. So the sun returned ten degrees, by which degrees it was gone down*'. (Isaiah 38:8, King James version, see also II Kings 20:11) Ahaz was king of Judah in around 730-700 BC. The story tells how the shadow on the dial miraculously moved backwards, as a sign from God that he would cure the sickness of Ahaz's son Hezekiah.

## *The Trypots*

The following quotation was mentioned in an excellent lecture by Tony Mann some years ago at BSHM meeting; it was what prompted me to read *Moby-Dick* and to write this article.

*Removing this hatch we expose the great trypots, two in number, and each of several barrels' capacity. When not in use, they are kept remarkably clean. Sometimes they are polished with soapstone and sand, till they shine within like silver punchbowls. During the night-watches some cynical old sailors will crawl into them and coil themselves away there for a nap. While employed in polishing them—one man in each pot, side by side—many confidential communications are carried on, over the iron lips. It is a place also for profound mathematical meditation. It was in the left hand try-pot of the Pequod, with the soapstone diligently circling round me, that I was first indirectly struck by the remarkable fact, that in geometry all bodies gliding along the cycloid, my soapstone for example, will descend from any point in precisely the same time.* [Ch.96, p.461]

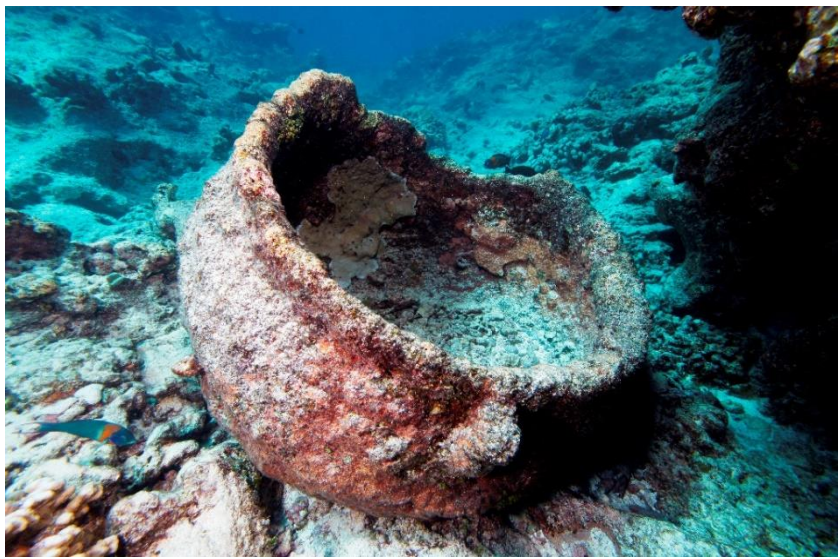

**Figure 3: Try-pot from the wreck of the whaling ship "Two Brothers", which sank in 1823 (public domain image created by US National Oceanic & Atmospheric Administration)**

The reference here is to the so-called tautochrone problem, to find the curve for which the time taken by a frictionless object sliding under gravity to the lowest point of the curve is independent of the starting point. A cycloid is the curve produced by a point on the circumference of a circle, or wheel, as it rolls along a straight line. If the circle has radius r, the line is the x-axis, and the point on the rim starts at the origin, then the cycloid consists of the points $r(\theta - \sin\theta, \theta - \cos\theta)$, where $\theta$ is the angle through which the circle has rotated. The first arch of the cycloid corresponds to the range $0 \leq \theta \leq 2\pi$ and it is to this shape that Melville is referring. As show in Figure 3, the curved shape of a trypot's cross-section, at least up to its widest point, could plausibly approximate a cycloid.

We do not know for sure who first discovered the cycloid, but a strong contender is Charles de Bovelles, who described the curve in his 1503 *Introductio in geometriam*. Bovelles was an influential French mathematician, whose *Géométrie en françoys* was the first scientific work in French ever printed. The tautochrone problem was first solved by Christiaan Huygens in 1659 – the required curve is a cycloid. The proof was given in his 1673 book *Horologium Oscillatorium*. If we take the inverted arch of the cycloid generated by a circle of radius *r*, then Huygens proved that the time of descent is $\pi\sqrt{\frac{r}{g}}$ (where *g* is acceleration due to gravity).

How did Melville learn about cycloids? They would not have been routinely taught in schools. However, though we do not know which edition of Euclid Melville had, cycloids were discussed in the introductions to several of the editions that could plausibly have belonged to him. Turpin mentions Potts' edition of 1845 as a possible contender (Turpin, 2015). Robert Potts (1805 – 1885) was an English mathematician; his edition of Euclid was widely used for much of the 19[th] century. However, I am not entirely persuaded by this: the introduction of Potts' edition (Potts, 1845) does mention cycloids, but only to assert (probably incorrectly) that they were invented by Galileo and that 'Pascal challenged the mathematicians of his day to prove some properties of the cycloid; but as the answers he received from Wallis and other eminent men were unsatisfactory, he himself gave the complete proofs of all the properties mentioned in his challenge'. Cycloids are not defined, nor are their properties, such as the tautochrone property, described. It is possible that the young Melville read this introduction, but without another source, this is not enough. As we shall see in the next section, Melville was very fortunate in one of his mathematics teachers; my conjecture is that he learned about cycloids from this teacher, as it is precisely the kind of thing one entertains bright pupils with.

## *Measurement*

The 19[th] century was a time of increased systematisation and formalisation across the sciences, with animal and plant species being ever more carefully classified, and taxonomies being standardised. (Mathematics too headed in this direction – perhaps we could call the quest of Russell and Whitehead's *Principia Mathematica* the *Moby-Dick* of mathematics, fatally undone by Gödel so soon after.) We will see later in our discussion of Maury's work (Section 6), that even such apparently quixotic phenomena as weather, ocean currents and the movements of whales, were starting to be analysed and pinned down in the manner of a Victorian lepidopterist affixing specimens in a glass case. In *Moby-Dick,* there are chapters on 'Cetology' in which Ishmael attempts a classification of whales (using the publishing terms Folio, Octavo and Duodecimo to divide the whales by size). The facetiousness of this classification speaks to Ishmael's understanding that, although we cannot hope to control nature and the sea by pure force of logic, we still, being human, need such structures to cling to, else we would drown in the ocean of information. Mathematics, in other words, is one way in which we navigate through the unknowable universe. Ishmael certainly values data – he engraves them upon his own body: '*The skeleton dimensions [of the whale] I*

*shall now proceed to set down are copied verbatim from my right arm, where I had them tattooed; as in my wild wanderings at that period, there was no other secure way of preserving such valuable statistics'.*[Ch.104, p.491]

Numbers and calculation play a key symbolic role in the novel. As long as Ahab is still rational enough to perform his daily navigational calculations, the voyage is not fatally imperilled. The inference is that, while mathematics and calculation can only give us partial understanding (we can get all the statistics on whale migrations that we like, but we can never predict what one whale will do), if we abandon them then that is truly a sign of madness. In the early part of the novel, the Captain makes regular calculations. In [Ch.34,p.161], Ahab, '*sitting in the lee quarter-boat, has just been taking an observation of the sun; and is now mutely reckoning the latitude on the smooth, medallion-shaped tablet, reserved for that daily purpose on the upper part of his ivory leg*'.

Several navigational instruments and calculations are mentioned, for example in Chapter 35, *The Mast-Head*, we hear of the '*binnacle deviations, azimuth compass observations and approximate errors*' spoken of while '*studying the mathematics aloft there in that bird's nest within three or four perches of the pole*'. Quadrants are used to calculate the latitude, of course, by determining the position of the sun. For longitude, a tool such as the 'patent chronometer' mentioned elsewhere in the novel, could be used. Alternatively, one could use lunar distances, computed with the new method of Bowditch described in his *Practical Navigator* (see Section 6).

Over the course of the ill-fated voyage, Ahab studies his maps and charts ever more obsessively in the hunt for the white whale, growing increasingly frustrated. Ahab's growing madness is shown when he rejects his calculations and instruments, beginning with the quadrant.

'*Then gazing at his quadrant, and handling, one after the other, its numerous cabalistical contrivances, he pondered again, and muttered: "Foolish toy! babies' plaything of haughty Admirals, and Commodores, and Captains; the world brags of thee, of thy cunning and might; but what after all canst thou do, but tell the poor, pitiful point, where thou thyself happenest to be on this wide planet, and the hand that holds thee: no! not one jot more! Thou canst not tell where one drop of water or one grain of sand will be to-morrow noon; and yet with thy impotence thou insultest the sun! Science! Curse thee, thou vain toy; and cursed be all the things that cast man's eyes aloft to that heaven, whose live vividness but scorches him, as these old eyes are even now scorched with thy light, O sun! Level by nature to this earth's horizon are the glances of man's eyes; not shot from the crown of his head, as if God had meant him to gaze on his firmament. Curse thee, thou quadrant!" dashing it to the deck, "no longer will I guide my earthly way by thee; the level ship's compass, and the level dead-reckoning, by log and by line; these shall conduct me, and show me my place on the sea. Aye," lighting from the boat to the deck, "thus I trample on thee, thou paltry thing that feebly pointest on high; thus I split and destroy thee!"'* [Ch.118,p.544]

The night that Ahab destroys the quadrant, there is a storm, which throws off the compasses. Seeing the '*crushed copper sight-tubes of the quadrant he had the day before dashed to the deck*', he then claims to be '*lord over the level loadstone yet*' [Ch.124, p.563] and makes his own. In the following chapter Ahab decides to use the log and line, which was even then a rarely-used means of measuring a ship's speed. When the line breaks too, Ahab is left with none of his navigational instruments – quadrant, compass, log and line are all lost. Logic and calculation have been rejected (evidence for Melville, who values these things, of Ahab's madness) and the *Pequod* is now sailing purely on the instinct of Ahab. Mathematics has been abandoned by the monomaniac, and we really are all at sea.

# 5. Melville's Education

At the time when Herman Melville was being educated, there were several different types of schools available in New York State; these are described in (Titus, 1980). Some private schools remained from colonial times (these were often called Latin schools). Then there were Academies, which came under the jurisdiction of the 'Board of Regents of the University of the State of New York'; they were essentially preparatory schools for College. Various pieces of legislation between 1780 and 1830 created (among other types) the so-called Union Free Schools, Select Schools and Lancastrian Schools. Legislation had not yet made school attendance compulsory, and only around half the children in the state had any formal schooling at all. Herman was enrolled at the New York Male High School in 1825, at the age of six. This was a Lancastrian School, which meant that instructors taught student 'monitors', who went on to teach other students. Pupils learned by rote in an extremely strict, punitive regime, in which the young Herman did not initially flourish. In August 1826 Melville's father described him in a letter as 'very backward in speech and somewhat slow in comprehension'. However, he found his feet, being selected as a monitor during the 1826-27 school year. By the time, at age ten, that he moved to the Columbia Grammar School, he was doing better. In a letter to his grandmother in October 1828, Melville writes that he is studying 'Geography, Grammar, Arithmetic, Writing, Speaking, Spelling and read [sic] in the Scientific class book' (Titus, 1980). Another letter from Allan at this time describes Herman as making 'more progress than formerly, and without being a bright Scholar, he maintains a respectable standing'. Allan also states that Herman has chosen Commerce as a 'favorite pursuit'. Parker believes that this is less about Herman actively choosing 'commerce', than 'assert[ing] some degree of control over his circumstances by solving problems in arithmetic'.

When the family moved to Albany in October 1830, Herman and his older brother Gansevoort were registered at the Albany Academy, recognised as the best school in the state, and one that had exceptionally good facilities – a well-equipped laboratory, an extensive library and state-of-the-art meteorological observation instruments (due to a law that required the thirty state academies to provide annual weather reports). Herman (now aged 11) studied geography, reading, spelling, penmanship, arithmetic, English grammar, and natural history (Parker, 1996). He was expected to learn by heart 'catechisms' in history, including that of ancient Greece and Rome, as well as classical biography and Jewish antiquities. In a report to the Board of Trustees of the academy in 1832, we are told that '*the afternoon is entirely devoted to Arithmetic. [E]ach student is employed one hour in his arithmetic lesson and engaged during the remainder of the Afternoon (1 hour 40 minutes) in entering sums into a large ciphering book, learning his lessons, or in examining the sums which are to be prepared for the succeeding day*' (Titus, 1980). Herman seems to have excelled in these lessons. In the Summer of 1831, he won a prize for being '*the first best in his class in ciphering books*'. Most biographers, though mentioning this, do so only in the context of pointing out that the award 'surprised everyone'. In fact, as Farmer notes (with some frustration!) one biographer 'immediately shifts from Melville's "mathematics examination" to the fact that the [poetry] book he received as a prize "ignited a spark of poetry in his soul" – arguably transposing Melville's math award into a sort of origin story for his literary future.' As a mathematician, my view is of course that mathematics itself is a form of poetry and that if Herman had a teacher who could show some of that poetry to his pupils (as seems likely, see below), then could it not have been just as much the mathematics as the prize that spurred his creativity?

Sadly, the family finances were in such a parlous state by October 1831, that Herman had to be withdrawn from school and find work. He had little education after this point, though in 1835 he joined the Albany Young Men's Association, a 'mutual improvement society', with debating clubs and access to many books for private study. Briefly in the Spring of 1835, Herman was enrolled in the Albany Classical School, and then in September 1836 he returned to the Albany Academy, where he enrolled in the Latin course. Again, he had to withdraw because of financial difficulties, and decided in 1837 to become a school teacher. He got a job at Sikes District School in Massachusetts, which lasted just one semester – this was one of several short-lived teaching positions. Finally, on November 12th, 1838 Melville enrolled at the Lansingburgh Academy where he completed 'what undoubtedly was a crash course in engineering and surveying during the two quarters he attended the school' (Titus, 1980). The object of these studies was seemingly to secure a job on the Erie Canal build. But, fortunately for literature, this job never materialised.

As so often, it may have been an inspirational teacher who gave Melville an interest in mathematics and science. The class in which he won the prize for ciphering was taught by Joseph Henry, Professor of Mathematics and Natural Philosophy at Albany Academy.

Joseph Henry was a well-known scientist who would go on to be the first Secretary of the Smithsonian. He was arguably the first build practical devices that used electromagnets – an early version of one of these was a bell that could be rung remotely via a battery connected with a wire. His experiments with electromagnetic relays were the basis of his telegraph machine, which he built at the Albany Academy – this laid the foundation for Morse's electrical telegraph. His discovery of inductance is the reason the modern SI unit of inductance is the henry (the discovery was independent of Faraday, though British readers will know that Faraday's work predates Henry's, which is why Faraday is credited with the discovery). In summary, Joseph Henry was a highly unusual school teacher.

Extant letters from Henry describe his difficulties in finding a good textbook, and a petition to the Academy's board (Farmer, 2016) 'to add a more advanced textbook for the Department's "higher students" just months before Melville won his ciphering award'. He was unsatisfied with his ultimate choice – Daboll's Arithmetic – dismissing it as 'a book of examples' 'good for nothing else'. Joseph Henry seems to have been a passionate and innovative teacher. He used visual aids like maps and charts, and demonstrated science experiments in class. He even used an exciting innovation, a blackboard, for teaching chemistry and mathematics. Henry's lectures on chemistry and physics, which were part of the advanced curriculum, were open to the public. He would often demonstrate his telegraph machine for awestruck pupils. Henry is the prime candidate, in my view, for the sort of teacher who may have taught his pupils about cycloids and the tautochrone problem.

Even after leaving school, Melville did not turn his back on technical subjects, as one may have expected. Farmer summarises his post-Academy years as follows: 'he worked at a bank, kept the books in a shop, and then taught mathematics before turning to a field that required him to learn geometry and trigonometry. Only after that failed – when he was unable to secure a position in engineering – Melville called upon "the only other talent he had" and submitted his first story' (Farmer, 2015). He did not train as a lawyer or go into business; he chose to train as a surveyor. This is the choice of someone aware that they are good at mathematics.

While we cannot know for sure which books Melville read, Turpin lists nine books 'touching on the subject of mathematics' that there is evidence for Melville having 'owned, borrowed or been familiar with' (Turpin, 2015). These include Euclid's *Elements*, *The American Practical Navigator*, by Nathaniel Bowditch, *Tracts on Mathematical and Philosophical Subjects*, by Charles Hutton, and *Daboll's Complete Schoolmaster's Assistant*, by Nathan Daboll. As we shall see in the next section, several of these make an appearance in *Moby-Dick*.

# 6. 'I have heard devils can be raised with Daboll's arithmetic'

In *Moby-Dick* four mathematicians, or at least names associated with mathematics, are mentioned: Euclid, Daboll, Bowditch and Maury. There are in fact mentions of others, but not in the context of mathematical work – Descartes is remarked upon in the context of philosophy and Pythagoras in the context of, let us say, good digestion (at sea, '*head winds are far more prevalent than winds from astern, that is, if you never violate the Pythagorean maxim*'[1]). In this section we discuss Daboll, Bowditch and Maury – anyone reading this will of course not need any introduction to Euclid!

## *Daboll*

Here is Stubb, the second mate of the Pequod, soliloquizing on the significance of the markings on a golden doubloon that Captain Ahab has offered as a reward for the first man to sight the white whale:

*Halloa! here's signs and wonders truly! [..] I'll get the almanack; and as I have heard devils can be raised with Daboll's arithmetic, I'll try my hand at raising a meaning out of these queer curvicues here with the Massachusetts calendar.*

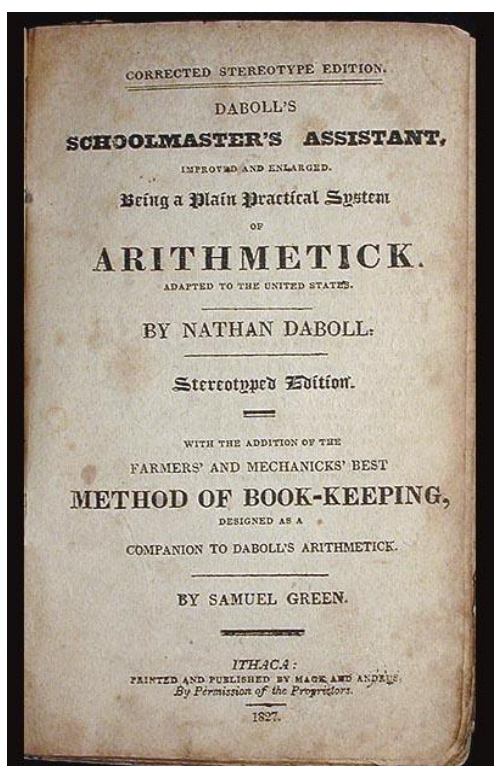

Figure 4: Front page of Daboll's Arithmetic (Nathan Daboll and Samuel Green, 1827 edition)

Generations of American schoolchildren would have been familiar with 'Daboll's arithmetic'. It was the most widely used text book in American schools from 1800 to at least 1850. The author, Nathan Daboll, was a mathematics teacher in Connecticut. His *Daboll's schoolmaster's assistant: being a plain, practical system of arithmetic, adapted to the United States*, first appeared in 1799, with an updated, expanded version published in 1814. Daboll went on to set up a navigation school, and his *Practical Navigator* was published (posthumously) in 1820. He also wrote, or contributed to, many almanacs. Florian Cajori, the American historian of mathematics, wrote that Daboll was one of the three 'great arithmeticians in America'.

Melville used Daboll's arithmetic as a schoolboy and almost certainly also as a teacher. The book is focused on practical skills that would be needed by merchants and tradesmen, but it contains some quite advanced material. The book starts with basic arithmetic, but contains such topics as annuities, pensions, compound interest, currency conversion, duodecimals, brokerage, extraction of roots, arithmetic and geometric progressions and even permutations. There are tables of units to help you work out the number of scruples in a dram, or the number of hogsheads in a pipe. The section on extraction of roots includes such questions as '*what is the biquadrate, or 4<sup>th</sup> root, of 19987173376*' and '*required the sursolid, or fifth root, of 6436343*'. Applications of extraction of roots are also given; for instance, given a cask of a certain shape with a bung of a certain diameter, a craftsman is ordered to make another of the same shape but to hold exactly twice as much; what will be the bung diameter and length of the new cask? For this calculation both cube roots and square roots are required. In the section on permutations we see questions like the following '*Seven gentlemen met at an inn, and were so well pleased with their host, and with each other, that they agreed to tarry so long as they, together with their host, could sit every day in a different position at dinner; how long must they have staid at said inn to have fulfilled their agreement?*' (The answer given is $110\frac{170}{365}$ years.) There is almost no geometry in the book; just a few pages about areas and volumes (for geometry we have Euclid, of course). The rules given are somewhat mysterious; for example to find the area of a circle, '*[m]ultiply half the diameter by half the circumference, and the product is the area; or, if the diameter is given without the circumference, multiply the square of the diameter by ,7854 and the product will be the area*'. No mention is made of π; to find the circumference of a circle given its diameter, the rule is as follows: '*As 7: is to 22:: so is the given diameter: to the circumference. Or, more exactly, as 113: is to 355:: &c. the diameter is found inversely*'. By this point one rather sympathises with Stubb – this 'arithmetic' does seem a little like alchemy.

## *Bowditch*

Here is what Ishmael has to say on Bowditch: '*Beware of enlisting in your vigilant fisheries any lad [...] who offers to ship with the Phædon instead of Bowditch in his head*'.

The Phædon, incidentally, likely refers to Plato's dialogue on the immortality of the soul, now more usually called Phædo, rather than the 1767 book by Moses Mendelssohn. In either case, the point Ishmael is making is that Nantucket ship-owners should avoid hiring philosophical daydreamers such as himself, if they can avoid it.

Here is Stubb again:

*'Book! you lie there; the fact is, you books must know your places. You'll do to give us the bare words and facts, but we come in to supply the thoughts. That's my small experience, so far as the Massachusetts calendar, and Bowditch's navigator, and Daboll's arithmetic go'.*

Nathaniel Bowditch (1773 – 1838) was an American mathematician. While studying bookkeeping as an apprentice, he taught himself algebra, calculus, Latin and French. This meant he was able to read many mathematical works; his translation of Laplace's *Méchanique Céleste* into English was an important contribution to the development of astronomy in America. However, by far his best-known work is *The American Practical Navigator*, first published in 1802 (Bowditch, 1802). Bowditch's Navigator does not just contain tables of numbers, and instruction on the use of navigational instruments. The 1802 edition has a large section on geometry, for example, with theorems and proofs (for example, the standard Euclid proof of Pythagoras's theorem is given, along with proofs of several circle theorems). It is a mathematical book totally different in spirit from, and hugely greater in scope than, Daboll's arithmetic.

This work arose from Bowditch's time at sea, first as a ship's clerk but later as part-owner. He had initially worked with the *Navigator* of John Hamilton Moore, but it had so many errors that he recalculated all the tables for himself, and ultimately decided to write his own book, *Bowditch's New American Practical Navigator*. Amazingly, the book has been in continuous publication since 1802 – regularly revised and updated, of course). The copyright was bought by the U.S. government in the 1860s, since which time it has published over 50 editions. In the Naval Observatory library in Washington, D.C. (to quote from the preface of the 2002 edition), 'one can find an original 1802 first edition of the *New American Practical Navigator*. One cannot hold this small, delicate, slipcovered book without being impressed by the nearly 200-year unbroken chain of publication that it has enjoyed. It sailed on U.S. merchantmen and Navy ships shortly after the quasi-war with France and during British impressment of merchant seamen that led to the War of 1812. It sailed on U.S. Naval vessels during operations against Mexico in the 1840's, on ships of both the Union and Confederate fleets during the Civil War, and with the U.S. Navy in Cuba in 1898. It went around the world with the Great White Fleet, across the North Atlantic to Europe during both World Wars, to Asia during the Korean and Vietnam Wars, and to the Middle East during Operation Desert Storm. It has circled the globe with countless thousands of merchant ships for 200 years. As navigational requirements and procedures have changed throughout the years, *Bowditch* has changed with them. Originally devoted almost exclusively to celestial navigation, it now also covers a host of modern topics. It is as practical today as it was when Nathaniel Bowditch, master of the *Putnam*, gathered the crew on deck and taught them the mathematics involved in calculating lunar distances.'

## *Maury*

In a chapter entitled 'The Chart', we see Captain Ahab, with furrowed brow, obsessively studying maps and charts. Melville writes about the theoretical possibility of using ships logs to construct migratory charts of the sperm whale. In fact precisely such a project was being undertaken while he was writing *Moby-Dick*. He added a footnote in production to that effect:

*Since the above was written, the statement is happily borne out by an official circular, issued by Lieutenant Maury, of the National Observatory, Washington, April 16th, 1851. By that circular, it appears that precisely such a chart is in course of completion; and portions of it are presented in the circular. "This chart divides the ocean into districts of five degrees of latitude by five degrees of longitude, perpendicularly through each of which districts are twelve columns for the twelve months; and horizontally through each of which districts are three lines; one to show the number of days that have been spent in each month in every district, and the two others to show the number of days on which whales, sperm or right, have been seen."*

Matthew Fontaine Maury (1806-1873) was an officer of the U.S. Navy, a pioneer of oceanography and meteorology, and an early proponent of what we now call Big Data. While serving at sea he became very interested in patterns of currents, winds and other nautical and meteorological phenomena. He became interested by the idea of producing charts showing these data, and when an injury forced him out of active service and into a desk job as head of the Depot of Charts and Instruments, he happened upon a way to create these charts. He discovered that in storage were log books from every American ship for a century, but nothing had been done with the vast amount of data they contained. He started to collect and collate these data, giving predictions for currents and winds at different times of the year in different places. They really are astonishing feats of calculation. In (Guarnieri, 2018) a wind diagram from Maury's *The Physical Geography of the Sea* is reproduced. Guarnieri explains that '[t]hese tables are the result of 1,159,353 distinct observations concerning the strength and direction of the wind, and, roughly, 100,000 observations on the barometric measurements at sea'. The charts he produced allowed navigators to plan new routes, which resulted in dramatic reductions in ocean-crossing times – the first voyage using Maury's data sailed a route that normally took 110 days in just 75 days.

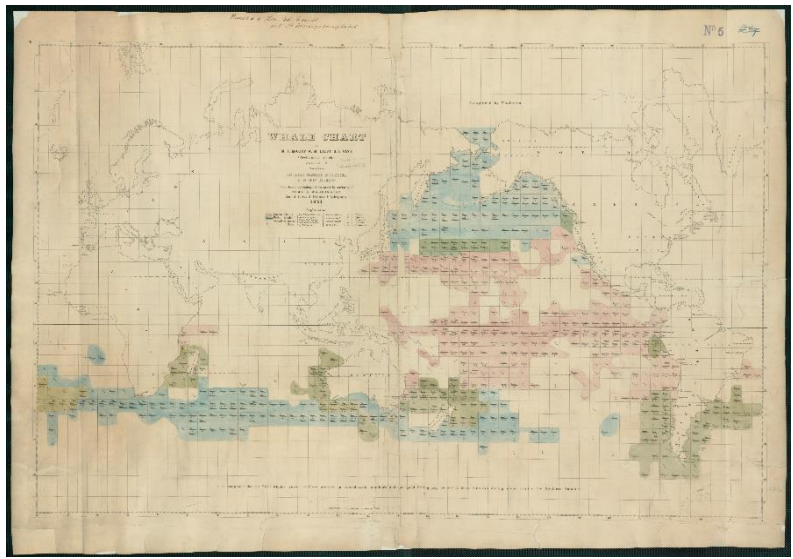

Figure 5: Matthew F. Maury. Whale Chart. Washington, D.C.: Naval Observatory, 1851. Map reproduction courtesy of the Norman B. Leventhal Map & Education Center at the Boston Public Library

The whale chart referred to in *Moby-Dick* is shown in Figure 5. It shows the distributions of several types of whales, indicating seasonal variations. It can be plainly seen from the chart that whales were mostly being found in the Pacific Ocean – this is unfortunately because whales had by that time been hunted almost to extinction in the Atlantic.

## *References by other authors*

Methods of schooling – the rote learning, the very small selection of texts (Daboll, Euclid) that were widely used for many decades – would have meant a common awareness of at least their names and some of the memory of having learnt from them, among many if not most of the educated American readers of the time. But this is not enough to explain their being mentioned by Melville; other authors did not do so. As a representative example, Nathaniel Hawthorne is one of the best-known of Melville's American literary contemporaries. Moreover, Melville and Hawthorne were close friends and neighbours, to the extent that *Moby-Dick* is dedicated to Hawthorne. However, a reader of Nathaniel Hawthorne's most famous novel, *The Scarlet Letter*, would struggle to find any mathematical allusions. In fact, a search of digitised copies of Hawthorne's novels shows that he uses the word 'mathematics', or 'mathematical', just twice in all his novels; the words 'Euclid' or 'Euclidean' do not appear at all (nor, unsurprisingly, do Daboll or Bowditch), while 'geometry' and 'arithmetic' each appear just once, essentially in the context of listing what a child has learnt – the 'rudiments of Latin and geometry', for example. Meanwhile in *Moby-Dick* alone, the words 'mathematics' or 'mathematical' are used six times, along with multiple uses of 'geometry', 'Euclid' or 'Euclidean', and 'arithmetic'. This does seem to be something characteristic of Melville rather than any great mathematical literacy in authors of the period. More importantly these are mathematically-minded uses, not as in, for instance (to cross the Atlantic briefly), the use of the word 'algebra' in Thackeray's Vanity Fair, published 1848, where he writes of one character that she 'no more comprehended sensibility than she did Algebra'.

# 7. Conclusion

Herman Melville, particularly in *Moby-Dick*, makes use of many more mathematical references and metaphors than other contemporary authors. We cannot explain this purely by proving that Melville knew more about Euclid than Hawthorne, for instance – though he almost certainly did. The choice of allusions used in a book is based on three things: the author's store of knowledge; their personal inclinations and interests; and what he or she wants us to infer about the narrator and about the larger themes of the novel. One would of course expect a book about a whaling voyage to feature more references to quadrants than, say, *Pride and Prejudice.* This, though, is not enough to explain the wealth of metaphors unconnected to navigation in the novel. We have seen that Melville had a good education in mathematics from an inspirational teacher. But transcending all this is what we hope this article has shown: that Melville actively enjoyed mathematics and mathematical ideas, and that this shines through in his work.